\documentclass[fleqn]{mat01}
\usepackage{times,amssymb,latexsym,amscd} 
\begin{document}

\setcounter{page}{391}
\firstpage{391}

\newtheorem{theore}{Theorem}
\renewcommand\thetheore{\arabic{theore}}
\newtheorem{theor}[theore]{\bf Theorem}
\newtheorem{lem}[theore]{Lemma}
\newtheorem{rem}[theore]{Remark}
\newtheorem{coro}[theore]{\rm COROLLARY}

\def\leq{\leqslant}
\def\geq{\geqslant}
\def\al{\alpha}
\def\be{\beta}
\def\ga{\gamma}
\def\C{\mathbb C}
\def\F{\mathbb F}
\def\N{\mathbb N}
\def\R{\mathbb R}
\def\Q{\mathbb Q}
\def\Z{\mathbb Z}
\def\eps{\varepsilon}
\def\ro{\varrho}

\newcommand{\Qbar}{{\overline{\Q}}}

\font\xxxx=mtsy at 10pt
\def\geq{\mbox{\xxxx{\,\char'025\,}}}
\def\leq{\mbox{\xxxx{\,\char'024\,}}}

\title{There are infinitely many limit points of the fractional
parts of powers}

\markboth{Art\= uras Dubickas}{Limit points of the fractional
parts of powers}

\author{ART\= URAS DUBICKAS}

\address{Department of Mathematics and Informatics, Vilnius
University, Naugarduko 24, Vilnius LT-03225, Lithuania\\
\noindent E-mail: arturas.dubickas@maf.vu.lt}

\volume{115}

\mon{November}

\parts{4}

\pubyear{2005}

\Date{MS received 20 May 2005; revised 19 July 2005}

\begin{abstract}
Suppose that $\al>1$ is an algebraic number and $\xi>0$ is a real
number. We prove that the sequence of fractional parts $\{\xi
\al^n\},$ $n =1,2,3,\dots,$ has infinitely many limit points
except when $\al$ is a PV-number and $\xi \in \Q(\al).$ For
$\xi=1$ and $\al$ being a rational non-integer number, this result
was proved by Vijayaraghavan.
\end{abstract}


\maketitle

\section{Introduction}

Let $\al>1$ and $\xi>0$ be real numbers. The problem of
distribution of the fractional parts $\{\xi \al^n\},$
$n=1,2,3,\dots,$ is a classical one. Some metrical results are
well-known. Firstly, for fixed $\al,$ the fractional parts $\{\xi
\al^n\},$ $n=1,2,3,\dots,$ are uniformly distributed in $[0,1)$
for almost all $\xi$ \cite{w}. Secondly, for fixed $\xi,$ the
fractional parts $\{\xi \al^n\},$ $n=1,2,3,\dots,$ are uniformly
distributed in $[0,1)$ for almost all $\al$ (see \cite{kok} and
also \cite{hl} for a weaker result). However, for fixed pairs
$\xi,\al,$ nearly nothing is known. Even the simple-looking
Mahler's question \cite{mah} about the fractional parts $\{\xi
(3/2)^n\},$ $n=1,2,3,\dots,$ is far from being solved. (See,
however, \cite{flp} and, for instance, see
\cite{ars,afs,bug,lond,dist} for more recent work on this
problem.)

One of the first results in this direction is due to
Vijayaraghavan, who proved that the set of limit points of the
sequence $\{(p/q)^n\},$ $n=1,2,3,\dots,$ where $p>q>1$ are
integers satisfying gcd$(p,q)=1,$ is infinite. In his note
\cite{vij} (see also \cite{vij1}) he gave two proofs of this fact:
one due to himself and another due to A~Weil. It was noticed later
that the questions of distribution of $\{\xi \al^n\},$
$n=1,2,3,\dots,$ for algebraic $\al$ are closely related to the
size of conjugates of $\al.$ The algebraic integers $\al>1$ whose
conjugates other than $\al$ itself are all strictly inside the
unit disc were named after Pisot and Vijayaraghavan and called
{\it PV-numbers} (see \cite{cas} and \cite{sal}).

The aim of this paper is to prove the following generalization of
the above mentioned result of Vijayaraghavan.

\begin{theor}[\!]
\label{pv} Let $\al>1$ be an algebraic number and let $\xi>0$ be a
real number. Then the set $\{\xi \al^n\},$ $n \in \N,$ has only
finitely many limit points if and only if $\al$ is a PV-number and
$\xi \in \Q(\al).$
\end{theor}

This theorem was already proved by Pisot in \cite{pis}. We give a
different proof by developing the method of Vijayaraghavan \cite{vij}.
In addition, we prove a stronger result for Salem numbers (see Lemma~3
below).

The `if' part of the theorem is well-known. Indeed, let
$\al=\al_1$ be a PV-number with conjugates, say, $\al_2, \dots,
\al_d.$ Assume that $\xi \in \Q(\al),$ that is, $\xi=(e_0+e_1 \al
+ \cdots + e_{d-1} \al^{d-1})/L$ with $e_0, \dots, e_{d-1} \in
\Z$ and $L \in \N.$ By considering the trace of $L\xi \al^n,$
namely, the sum over its conjugates, we have
\begin{align*}
{\rm Tr}(L\xi \al^n) &= e_0 {\rm Tr}(\al^n)+\cdots+e_{d-1} {\rm
Tr}(\al^{n+d-1})\\[.2pc]
&=L[\xi \al^n]+L\{\xi \al^n\}+e_0 \sum_{j=2}^d \al_j^n + \cdots +
e_{d-1} \sum_{j=2}^d \al_j^{n+d-1}.
\end{align*}
Since ${\rm Tr}(L\xi \al^n)-L[\xi \al^n]$ is an integer and, for
each fixed $k,$ the sum $\sum_{j=2}^d \al_j^{n+k}$ tends to zero
as $n \to \infty,$ we deduce that the set of limit points of
$\{\xi \al^n\},$ $n \in \N,$ is a subset of $\{0,1/L,\dots,
(L-1)/L, 1\}.$

So in the proof below we only need to prove the `only if' part,
namely, that in all other cases the set of limit points of $\{\xi
\al^n\},$ $n \in \N,$ is infinite.

We remark that the theorem does not apply to transcendental
numbers $\al>1.$ It is not known, for instance, whether the sets
$\{e^n\},$ $n \in \N,$ and $\{\pi^n\},$ $n \in \N,$ have one or
more than one limit point.

In some cases the theorem cannot be strengthened. Suppose, for
instance, that $\al$ is a rational integer $\al=b \geq 2$ (which
is a PV-number) and $\xi=\sum_{k=0}^{\infty}b^{-k!}$ (which is a
transcendental Liouville number, so $\xi \notin\Q(b)=\Q$). Then
the set of limit points of the sequence $\{\xi b^n\},$
$n=1,2,3,\dots,$ is $\{0,b^{-1},b^{-2},b^{-3},\dots\}.$
Evidently, this set is countable.

\section{Sketch of the proof and auxiliary results}

From now on, let us assume that $\al=\al_1>1$ is a fixed algebraic
number with conjugates $\al_2, \dots, \al_d$ and with minimal
polynomial $a_d z^d+a_{d-1} z^{d-1} + \cdots + a_0 \in \Z[z].$ Set
$L(\al)=|a_0|+|a_1|+\cdots+|a_d|.$ Suppose that $\xi>0$ is a real
number satisfying $\xi \notin \Q(\al)$ in case $\al$ is a
PV-number.

Recall that an algebraic integer $\al>1$ is called a {\it Salem
number} if its conjugates are all in the unit disc $|z| \leq 1$
with at least one conjugate lying on $|z|=1$. The next lemma is
part of Theorem~1 in \cite{dist}. (Here and below, $\|x\|:=\min
(\{x\}, 1-\{x\}).$)

\begin{lem}
\label{pvs} Let $\al>1$ be a real algebraic number and let $\xi>0$
be a real number. If $\|\xi \al^n\|<1/L(\al)$ for every $n \in \N$
then $\al$ is a PV-number or a Salem number and $\xi \in \Q(\al).$
\end{lem}

Suppose that the set $S$ of limit points of $\{\xi \al^n\},$ $n
\in \N,$ is finite, say, $S=\{\mu_1,\mu_2,\dots,\mu_q\}.$ With
this assumption, we will show in \S3 that, for any $\eps>0,$ there
exist three positive integers $m, r, L,$ where $m>r,$ such that
\begin{equation*}
\|L\xi(\al^m - \al^r)\al^n\| < 2\eps
\end{equation*}
for every $n \in \N.$ Taking $\eps<1/2L(\al),$ by Lemma~\ref{pvs},
we conclude that $\al$ is a PV-number or a Salem number and
$L\xi(\al^m-\al^r) \in \Q(\al),$ that is, $\xi \in \Q(\al).$

However, the case when $\al$ is a PV-number and $\xi \in \Q(\al)$
is already treated in the `if' part of the theorem. So the only
case that remains to be settled is when $\al$ is a Salem number
and $\xi \in \Q(\al).$ We will then prove even more than required.

\begin{lem}
\label{sal} Suppose that $\al$ is a Salem number and $\xi>0$
belongs to $\Q(\al).$ Then there is an interval $I=I(\xi,\al)
\subset [0,1]$ of positive length such that each point $\zeta \in
I$ is a limit point of the set $\{\xi \al^n\},$ $n \in \N.$
\end{lem}

We will prove Lemma~\ref{sal} in \S4. Finally, recall that the
sequence $b_1, b_2, b_3,\dots$ is called {\it ultimately
periodic} if there is a $t \in \N$ such that $b_{n+t}=b_n$ for
every $n \geq n_0$. If $n_0=1,$ then the sequence $b_1, b_2,
b_3,\dots$ is called {\it purely periodic}. The next lemma was
proved in \cite{pal}. It will be used in the proof of
Lemma~\ref{sal}.

\begin{lem}
\label{puper} Let $d, L \in \N$ and $A_{d-1}, \dots, A_0 \in \Z,$
where $A_0 \ne 0$ and $\gcd(A_0, L)=1.$ Then the sequence of
integers $b_1, b_2, b_3, \dots$ satisfying the linear recurrence
sequence
\begin{equation*}
b_{k+d}+A_{d-1}b_{k+d-1}+\cdots+A_1b_{k+1}+A_0b_k=0,
\end{equation*}
where $k=1,2,3,\dots,$ is purely periodic modulo $L.$
\end{lem}

\section{Differences of fractional parts are close to an integer}

Suppose that the set $S$ of limit points of $\{\xi \al^n\},$ $n
\in \N,$ is finite, say, $S=\{\mu_1, \dots, \mu_q\}.$ Let
$\mathcal {D}(S)$ be the set of all differences $\mu_i-\mu_j,$
where $\mu_i \geq \mu_j$ belong to $S.$ It is possible that the
set $S \cup \mathcal{D}(S)$ (which is a subset of $[0,1]$)
contains some rational numbers. For instance, $\mathcal{D}(S)$
always contains $0.$ Let $L$ be the least common multiple of the
denominators of all rational numbers that belong to $S \cup
\mathcal{D}(S).$ (Of course, $L:=1$ if $(S \cup \mathcal{D}(S))
\cap \Q =\{0\}$ or $\{0,1\}.$)

Consider the set $S_L$ of limit points of $\{L\xi \al^n\},$ $n \in
\N.$

\begin{lem}
\label{lll} $S_L$ is a subset of $\{0, \{L\mu_1\}, \dots,
\{L\mu_q\}, 1\}$.
\end{lem}

\begin{proof} Note that
\begin{equation*}
L\{\xi \al^{n}\} - \{L \xi \al^{n}\}=[L\xi \al^{n}] - L [\xi
\al^{n}]
\end{equation*}
is a non-negative integer. Therefore each element of $S_L$ is of
the form $L \mu_i - n_i$ with integer $n_i \geq 0.$ Evidently,
$S_L-\{0,1\}$ is a subset of the interval $(0,1).$ Consequently,
$n_i=[L \mu_i]$ for each $\mu_i$ satisfying $L\mu_i \notin \Z.$
This proves the lemma.\hfill $\Box$
\end{proof}

\begin{lem}
\label{dif} The set $(S_L \cup \mathcal{D}(S_L)) \cap \Q$ is
either $\{0, 1\}$ or $\{0\}.$
\end{lem}

\begin{proof}
Of course, for any rational $\mu_i,$ by the definition of $L,$ we
have $\{L \mu_i\} =0.$ By Lemma~\ref{lll}, we deduce that $S_L
\subset \{0, \dots,\{L\mu\},\dots ,1\},$ where $\mu$ runs over
every irrational element of $S,$ so that $S_L \cap \Q \subset
\{0,1\}.$ The difference $\{L \mu_i\}-\{L\mu_j\}=
L(\mu_i-\mu_j)-[L\mu_i]+[L\mu_j],$ where $\mu_i, \mu_j \in S,$ is
either irrational or, by the definition of $L,$ an integer. Hence,
$\mathcal{D}(S_L)$ contains at most two rational elements, namely,
$0$ and $1.$ This proves the lemma.\hfill $\Box$
\end{proof}

Write
\begin{equation*}
x_n = [L\xi \al^n] \quad {\rm and} \quad y_n= \{L\xi \al^n\}.
\end{equation*}
Then, as $a_0 \al^n + a_1 \al^{n+1}+\cdots + a_d \al^{n+d}=0,$ we
set
\begin{align*}
s_n:=a_0 y_n + a_1 y_{n+1} + \cdots + a_d y_{n+d} = -a_0 x_n-a_1
x_{n+1} - \cdots - a_d x_{n+d}.
\end{align*}
So $s_n$ belongs to a finite set of integers for each $n \in \N.$
(We remark that a key result which was proved in \cite{lond} is
that the sequence $s_1, s_2, s_3, \dots$ is not ultimately
periodic, unless $\al$ is a PV-number or a Salem number and $\xi
\in \Q(\al).$ Lemma~\ref{pvs} given in \S2 is an easy consequence
of this result.)

Suppose that $S_L$ contains $g$ irrational elements. We denote
$S_L^*=S_L-\{0,1\}.$ By Lemma~\ref{dif}, the set $S_L$ contains at
most two rational elements $0$ and $1.$ Hence $S_L$ contains at
most $g+2$ elements. By Lemma~\ref{dif} again, the numbers
$\eta-\eta',$ where $\eta, \eta' \in S_L,$ $\eta>\eta',$ are all
irrational except (possibly) when $(\eta, \eta')=(1,0).$

Set
\begin{equation*}
\tau = \min \|a_d (\eta-\eta')\|,
\end{equation*}
where the minimum is taken over every pair $\eta, \eta' \in S_L^*
\cup \{0,1\},$ where $\eta>\eta',$ except for the pair $(\eta,
\eta')=(1,0).$ Since all these differences are irrational, we have
$0<\tau<1/2.$

Recall that $s_n=a_0 y_n + \cdots + a_d y_{n+d}$ is an integer,
where $y_n=\{L \xi \al^n\}.$ Fix $\eps$ in the interval
$0<\eps<\tau/2L(\al)<1/4L(\al).$ Then the intervals $(\eta-\eps,
\eta+\eps),$ where $\eta \in S_L^* \cup \{0,1\},$ are disjoint.
Furthermore, there is an integer $N$ so large that $y_n$ lies in
an $\eps$-neighbourhood of $\eta=\eta_n \in S_L$ for each $n \geq
N.$ We will write $\eta_n$ for the element of $S_L$ closest to
$y_n.$

Consider the vectors $Z_h:=(\eta_h, \eta_{h+1}, \dots,
\eta_{h+d})$ for $h=N, N+1, \dots.$ There are at most
$(g+2)^{d+1}$ different vectors in $S_L^{d+1}.$ So there are two
integers, say, $m$ and $r$ satisfying $m>r \geq N,$ such that
$Z_m=Z_r.$ Subtracting $s_{r+n}$ from $s_{m+n}$ yields
\begin{align*}
s_{m+n}-s_{r+n} &= a_{0}(y_{m+n}-y_{r+n})+ \cdots +
a_{d-1}(y_{m+d-1+n}-y_{r+d-1+n})\\[.2pc]
&\quad\, +a_d(y_{m+d+n}-y_{r+d+n})
\end{align*}
for $n=0,1,\dots.$ Writing $y_h=\eta_h+(y_h-\eta_h)$ and using
$|y_h-\eta_h|<\eps,$ we deduce that
\begin{align*}
&\|a_0(\eta_{m+n}-\eta_{r+n})+\cdots+a_{d-1}(\eta_{m+d-1+n}-
\eta_{r+d-1+n})\\[.2pc]
&\quad\, + a_{d}(\eta_{m+d+n}-\eta_{r+d+n})\|<2 \eps L(\al) <
\tau.
\end{align*}

We next claim that the difference $\eta_{m+n}-\eta_{r+n}$ belongs
to the set $\{0,1,-1\}$ for each $n \geq 0.$ Since $Z_m=Z_r,$ we
have $\eta_{m+n}=\eta_{r+n}$ for every $n=0,1,\dots, d.$ For the
contradiction, assume that $l$ is the smallest positive integer
for which $\eta_{m+d+l} - \eta_{r+d+l} \notin \{0,1,-1\}.$ In
particular, this implies that $\eta_{m+j+l}-\eta_{r+j+l} \in \Z$
for $j=0,1,\dots,d-1.$ Hence
\begin{align*}
&\|a_d(\eta_{m+d+l}-\eta_{r+d+l})\|\\[.2pc]
&\quad\, = \|a_0(\eta_{m+l}- \eta_{r+l})+\cdots+
a_{d}(\eta_{m+d+l}- \eta_{r+d+l})\| < \tau.
\end{align*}
By the choice of $\tau,$ this is impossible, unless
$\eta_{m+d+l}=\eta_{r+d+l}$ or $\{\eta_{m+d+l},
\eta_{r+d+l}\}=\{0,1\}.$ However, in both cases, we have
$\eta_{m+d+l} - \eta_{r+d+l} \in \{0,1,-1\},$ a contradiction.

Note that, since $\eta_h \in [0,1],$ we have
$\eta_{m+n}-\eta_{r+n} \in \{0,1,-1\}$ if and only if either
$\eta_{m+n}=\eta_{r+n}$ or $\{\eta_{m+n}, \eta_{r+n}\}=\{0,1\}.$
Obviously, $\eta_{m+n}=\eta_{r+n}$ implies that the difference
between the fractional parts $y_{m+n}=\{L\xi \al^{m+n}\}$ and
$y_{r+n}=\{L\xi \al^{r+n}\}$ is smaller than $2\eps.$ The
alternative case, namely, $\{\eta_{m+n}, \eta_{r+n}\}=\{0,1\}$
occurs when one of the numbers $\{L\xi \al^{m+n}\}, \{L\xi
\al^{r+n}\}$ is smaller than $\eps$ and another is greater than
$1-\eps.$ So, in both cases, we have
\begin{equation*}
\|L\xi \al^{m+n}-L\xi \al^{r+n}\| < 2\eps
\end{equation*}
for every $n \geq 0.$ Thus, we established the existence of three
positive integers $m, r, L,$ where $m>r,$ such that $\|L\xi
(\al^{m}-\al^{r})\al^n\| < 2\eps$ for every $n \in \N$ (as
required in \S2).

\section{Salem numbers}

In this section, we will prove Lemma~\ref{sal} and thus complete
the proof of the theorem. Suppose that $\al$ is a Salem number.
Let us write the conjugates of $\al$ in the form $\al^{-1},$ ${\rm
e}^{\phi_1 \sqrt{-1}}, {\rm e}^{-\phi_1 \sqrt{-1}}, \dots, {\rm
e}^{\phi_m \sqrt{-1}}, {\rm e}^{-\phi_m \sqrt{-1}},$ where
$d=2m+2$ and where the arguments $\phi_1, \dots, \phi_{m}$ belong
to the interval $(0,\pi).$ As above, we set $\xi=(e_0+e_1 \al +
\cdots + e_{d-1} \al^{d-1})/L$ with $e_0, \dots, e_{d-1} \in \Z$
and $L \in \N.$ Now, by considering the trace of $L\xi \al^n,$ we
have
\begin{align*}
{\rm Tr}(L\xi \al^n) &= e_0 {\rm Tr}(\al^n)+\cdots+e_{d-1} {\rm
Tr}(\al^{n+d-1})\\[.2pc]
&= L[\xi \al^n]+L\{\xi \al^n\}+ (e_0+\cdots+e_{d-1} \al^{-d+1})\al^{-n}\\[.2pc]
&\quad\, + 2 \sum_{j=1}^m (e_0 \cos(n\phi_j)+e_1\cos((n+1)
\phi_j)\\[.2pc]
&\quad\, +\cdots +e_{d-1}\cos((n+d-1)\phi_j)).
\end{align*}

Setting
\begin{equation*}
U(z):=e_0+e_1\cos z +\cdots +e_{d-1} \cos ((d-1)z)
\end{equation*}
and
\begin{equation*}
V(z):=e_1 \sin z + \cdots +e_{d-1} \sin((d-1)z),
\end{equation*}
we can write
\begin{align*}
&\sum_{j=1}^m (e_0 \cos(n\phi_j)+e_1\cos((n+1)\phi_j)+\cdots+
e_{d-1}\cos((n+d-1)\phi_j))\\[.2pc]
&\quad\, = \sum_{j=1}^m (U(\phi_j) \cos (n\phi_j)
-V(\phi_j)\sin(n\phi_j)).
\end{align*}
It follows that the sum
\begin{align*}
&L\{\xi \al^n \}+2\sum_{j=1}^m (U(\phi_j) \cos (n\phi_j)
-V(\phi_j)\sin(n\phi_j))\\[.2pc]
&\quad\, = {\rm Tr}(L\xi
\al^n)-L[\xi\al^n]-(e_0+\cdots+e_{d-1}\al^{-d+1})\al^{-n}
\end{align*}
is close to an integer for each sufficiently large $n.$ Moreover,
the sequence of integers $b_n:={\rm Tr}(L\xi \al^n),$
$n=1,2,3,\dots,$ satisfies the linear recurrence sequence
\begin{equation*}
a_d b_{n+d}+a_{d-1}b_{n+d-1}+\cdots+a_0b_n=0
\end{equation*}
for $n=1,2,3,\dots.$ The fact that $\al$ is a Salem number
implies that $a_d=a_0=1,$ so we can apply Lemma~\ref{puper}. It
follows that the sequence $b_n={\rm Tr}(L\xi \al^n),$
$n=1,2,\dots,$ is purely periodic modulo $L.$ Let $q$ be the
length of its period, so that the numbers $b_q, b_{2q},
b_{3q},\dots$ are all equal modulo $L.$ Then, there is an integer
$\ell$ in the range $0\leq \ell \leq L-1$ such that
\begin{equation*}
L\{\xi \al^{qn}\}+2R_n \to \ell
\end{equation*}
as $n \to \infty.$ Here, $R_n:=\sum_{j=1}^m (U(\phi_j) \cos
(qn\phi_j) -V(\phi_j)\sin(qn\phi_j)).$\pagebreak

Note that
\begin{align*}
&U(\phi) \cos (qn \phi) -V(\phi)\sin(qn\phi)\\[.2pc]
&\quad\, = \sqrt{U(\phi)^2+ V(\phi)^2} \cos(qn\phi -
\arctan(V(\phi)/U(\phi)))
\end{align*}
for each real number $\phi.$ The numbers $\phi_1, \dots,
\phi_{m}$ and $\pi$ are linearly independent over $\Q$ (see, for
example, p.~32 of \cite{sal}). Hence, by Kronecker's theorem, for
arbitrary $m$ numbers $\theta_1, \dots, \theta_m \in [-1,1]$
there is an $n \in \N$ such that the value of $U(\phi_j) \cos
(qn\phi_j) -V(\phi_j)\sin(qn\phi_j)$ lies close to $\theta_j
\sqrt{U(\phi_j)^2+V(\phi_j)^2}$ for every $j=1,\dots,m.$

It follows that the sequence $R_n,$ $n=1,2,\dots,$ is dense in
the interval $[-H,H],$ where $H=\sum_{j=1}^m \sqrt{U(\phi_j)^2+
V(\phi_j)^2}.$ Clearly, $H>0,$ because $H=0$ yields $e_0+e_1 \al'
+ \cdots + e_{d-1} \al^{\prime d-1}=U(\phi_1)+\sqrt{-1}
V(\phi_1)=0,$ where $\al'={\rm e}^{\phi_1 \sqrt{-1}}$ is conjugate
to $\al.$ This is impossible, because the degree of $\al'$ over
$\Q$ equals $d.$

Now, since $\{\xi\al^{qn}\}+2R_n/L \to \ell/L$ as $n \to \infty$
and since the values of $2R_n,$ $n \in \N,$ are dense in
$[-2H,0],$ we see that there exists an interval $[\ell/L,
\ell/L+\delta],$ where $\delta$ is a positive number, such that
every $\zeta \in [\ell/L,\ell/L+\delta]$ is a limit point of the
set $\{\xi \al^{qn}\},$ $n \in \N.$ This completes the proof of
Lemma~\ref{sal}.\hfill $\Box$\vspace{.6pc}

See also \cite{arch,za1,za2} for other recent results concerning
integer and fractional parts of Salem numbers.

\section*{Acknowledgement}

This research was supported in part by the Lithuanian State
Studies and Science Foundation.

\end{document}